\documentclass[11pt]{article}
\usepackage{amsmath}
\usepackage{dsfont}
\usepackage{mathrsfs}
\usepackage{amsmath,amssymb}
\usepackage{amsfonts}
\usepackage{hyperref}
\usepackage{amsthm}
\usepackage{graphicx}
\usepackage{subfigure}
\usepackage{xcolor}
\usepackage{appendix}
\usepackage{cite}

\usepackage{tikz}

\usetikzlibrary{calc,shapes,arrows}

\usepackage{pgfplots}

\usetikzlibrary{patterns}

\usepgfplotslibrary{dateplot}

\renewcommand{\qed}{\hfill\small{$\square$}\normalsize}

\hfuzz=\maxdimen
\theoremstyle{definition}
\newtheorem{lemma}{Lemma}[section]
\newtheorem{definition}[lemma]{Definition}
\newtheorem{proposition}[lemma]{Proposition}
\newtheorem{theorem}[lemma]{Theorem}
\newtheorem{corollary}[lemma]{Corollary}

\newtheorem{remark}{Remark}

\newtheorem{conjecture}{Conjecture}
\numberwithin{equation}{section}
\renewcommand{\proof}{\textbf{Proof. }}
\renewcommand{\qed}{\hfill\small{$\square$}\normalsize}

\DeclareFixedFont{\Acknowledgment}{OT1}{cmr}{bx}{n}{14pt}
\begin{document}

\title{\bf Thurston's sphere packings on 3-dimensional manifolds, I}
\author{Xiaokai He, Xu Xu}
\maketitle

\begin{abstract}
Thurston's sphere packing on a 3-dimensional manifold is an analogue of Thurston's circle packing on a surface,
the rigidity of which has been open for many years.
In this paper,
we prove that Thurston's Euclidean sphere packing is locally determined by combinatorial scalar curvature up to scaling,
which generalizes Cooper-Rivin-Glickenstein's local rigidity for tangential sphere packing on 3-dimensional manifolds.
We also prove the infinitesimal rigidity that Thurston's Euclidean sphere packing can not be deformed (except by scaling) while keeping the combinatorial Ricci curvature fixed. The main tool used in the proof is a variational principle with the discrete Hilbert-Einstein functional as the action functional.
\end{abstract}



\section{Introduction}
In the study of hyperbolic metrics on 3-dimensional manifolds,
Thurston (\cite{T1}, Chapter 13) proved the Koebe-Andreev-Thurston theorem for circle packings
with non-obtuse intersection angles on surfaces,
the rigidity part of which states that a circle packing is determined by discrete curvatures on a triangulated surface.
Thurston's work generalizes Andreev's work on the sphere \cite{Andr1, Andr2}
and Koebe's work \cite{K1} for tangential circle packing on the sphere.
For a proof of Koebe-Andreev-Thurston theorem, see \cite{BSp,CL, DV, MR, St, T1, Z}.

A longstanding problem is whether the circle packings on surfaces could be generalized to three or higher dimensional manifolds
with the properties of circle packings on surfaces preserved, including the rigidity and existence.
Cooper-Rivin \cite{CR} introduced
the tangential sphere packing on a 3-dimensional manifold,
which was proved to be locally rigid \cite{CR,G1,G2,R} and globally rigid \cite{X1}.
Glickenstein \cite{G3} and Thomas \cite{T} studied the general discrete conformal deformations of $3$-manifolds and
obtained some rigidity results in some special cases.
The deformations of tangential sphere packing metrics by discrete curvature flows were studied in \cite{GH,GJS,GM,GX1,GX2,GX4,G1,G2}.
Thurston's sphere packing on a $3$-dimensional manifold is
an analogue of Thusrton's circle packing on a surface
and a natural generalization of Cooper-Rivin's tangential sphere packing on a $3$-dimensional manifold,
the rigidity of which has been open for many years.
In this paper, we study the infinitesimal and local rigidity of Thurston's Euclidean sphere packings on 3-dimensional manifolds.


Suppose $M$ is a 3-dimensional connected closed manifold with a triangulation $\mathcal{T}=\{V,E,F,T\}$,
where the symbols $V,E,F,T$
represent the sets of vertices, edges, faces and tetrahedra, respectively.
We denote a triangulated 3-manifold as $(M, \mathcal{T})$
and denote the simplices in $\mathcal{T}$ as $\{ij\cdots k\}$, where $i,j,k$ are vertices in $V$.

\begin{definition}\label{definition of Thurston's sphere packing metric}
Suppose $(M, \mathcal{T})$ is a triangulated closed connected 3-dimensional manifold with a weight $\Phi: E\rightarrow [0, \frac{\pi}{2}]$.
Thurston's Euclidean sphere packing metric is a map $r:V\rightarrow (0,+\infty)$ such that
\begin{description}
  \item[(1)] The length of the edge $\{ij\}\in E$ between the
vertices $i$ and $j$ is
\begin{equation}\label{Thurston Euclidean length introduction}
l_{ij}=\sqrt{r_{i}^2+r_{j}^2+2r_ir_j\cos \Phi_{ij}}.
\end{equation}
  \item[(2)] The lengths $l_{ij},l_{ik},l_{il},l_{jk},l_{jl},l_{kl}$ determine a nondegenerate Euclidean
tetrahedron for each tetrahedron $\{ijkl\}\in T$.
\end{description}
\end{definition}

\begin{remark}
Thurston's sphere packing metric can also be defined for the hyperbolic background geometry with the length of the edge $\{ij\}\in E$ replaced by
\begin{equation*}
l_{ij}=\cosh^{-1}(\cosh r_{i}\cosh r_j+\sinh r_i\sinh r_j\cos \Phi_{ij})
\end{equation*}
and the lengths $l_{ij},l_{ik},l_{il},l_{jk},l_{jl},l_{kl}$ determining a nondegenerate hyperbolic tetrahedron.
\end{remark}

A triangulated 3-manifold $(M, \mathcal{T})$ with a weight $\Phi: E\rightarrow [0, \frac{\pi}{2}]$ is denoted as $(M, \mathcal{T}, \Phi)$ in the following.
When $\Phi\equiv 0$, Thurston's sphere packing metric is
Cooper-Rivin's tangential sphere packing metric \cite{CR}. In this paper, we focus on Thurston's sphere packing metrics in Definition \ref{definition of Thurston's sphere packing metric}. For simplicity, we use Thurston's sphere packing metrics
to denote Thurston's Euclidean sphere packing metrics in the following, if it causes no confusion in the context.
The geometric meaning of Thurston's sphere packing metric is as follows.
Given a Thurston's sphere packing metric $r$ on $(M, \mathcal{T}, \Phi)$,
a sphere $S_i$ of radius $r_i$ is attached to each vertex $i\in V$, with the intersection angle of $S_i, S_j$ attached to
the end points $i, j$ of an edge $\{ij\}$ given by $\Phi_{ij}$.
For each topological tetrahedron in $T$, the sphere packing metric determines a nondegenerate Euclidean tetrahedron by Definition \ref{definition of Thurston's sphere packing metric}.
Isometrically gluing the Euclidean tetrahedra
along the faces gives rise to a piecewise linear $3$-manifold, which has singularities along the edges and vertices.
Here we use  singularities to denote the points where the piecewise linear metric is not smooth. 
The combinatorial Ricci and scalar curvature are used to describe the singularities along the edges and vertices respectively.



\begin{definition}[\cite{Re}]\label{defn comb Ricci curv}
Suppose $(M, \mathcal{T})$ is a triangulated 3-manifold with a piecewise linear metric.
The combinatorial Ricci curvature along the edge $\{ij\}\in E$ is defined to be
\begin{equation*}
\begin{aligned}
K_{ij}=2\pi-\sum_{\{ijkl\}\in T}\beta_{ij,kl},
\end{aligned}
\end{equation*}
where $\beta_{ij,kl}$ is the dihedral angle along the edge $\{ij\}$ in the tetrahedron $\{ijkl\}$ and the summation is taken over
tetrahedra with $\{ij\}$ as a common edge.
\end{definition}

For a Riemannian manifold  $(M, g)$, Cheeger-M\"{u}ller-Schrader \cite{CMS} proved that $K_{ij} l_{ij}$ converges in measure to the scalar curvature measure $RdV$,
where $R$ is the scalar curvature and $dV$ is the volume form of $(M, g)$.

We have the following infinitesimal rigidity for Thurston's sphere packing
with respect to the combinatorial Ricci curvature in Definition \ref{defn comb Ricci curv}.

\begin{theorem}\label{main theorem infinitesimal rigidity}
Suppose $(M, \mathcal{T})$ is a triangulated closed connected 3-manifold
with a weight $\Phi: E\rightarrow [0,\frac{\pi}{2}]$ satisfying
\begin{equation}\label{condition 1}
\begin{aligned}
\Phi_{ij}+\Phi_{ik}+\Phi_{jk}> \pi, \forall \{ijk\}\in F,
\end{aligned}
\end{equation}
or
\begin{equation}\label{condition 2}
\begin{aligned}
\Phi\equiv C\in [0,\frac{\pi}{2}].
\end{aligned}
\end{equation}
Then Thurston's sphere packings on $(M, \mathcal{T}, \Phi)$
can not be deformed (except by scaling) while keeping the combinatorial Ricci curvatures along the edges fixed.
\end{theorem}


Using the combinatorial Ricci curvature, one can further define the combinatorial scalar curvature for Thurston's sphere packing metrics.
\begin{definition}[\cite{G3}]\label{def of scalar curvature}
Suppose $(M, \mathcal{T}, \Phi)$ is a weighted triangulated connected closed 3-manifold with a nondegenerate sphere packing metric
$r:V\rightarrow (0,+\infty)$.
The combinatorial scalar curvature at a vertex $i$ for the sphere packing metric $r$ is defined to be
\begin{equation*}
K_i=\sum_{j\sim i}K_{ij}\cos\tau_{ij},\\
\end{equation*}
where $\tau_{ij}$ is the inner angle facing the edge with length $r_j$ in the triangle formed by three edges with lengths $r_i,r_j,l_{ij}$.
\end{definition}


\begin{remark}
The combinatorial scalar curvature in Definition \ref{def of scalar curvature}
comes from the first variation of the Einstein-Hilbert functional with respect to Thurston's sphere packing metrics.
Please refer to Section \ref{section 4} for more details.
Glickenstein \cite{G3} once defined the combinatorial scalar curvature for general discrete conformal metrics on 3-manifolds.
The combinatorial scalar curvature $K_i$ in Definition \ref{def of scalar curvature} is essentially the scalar curvature of Glickenstein divided by $r_i$  in the case of Thurston's sphere packing metrics.
If $\Phi_{ij}\equiv0$ for any $\{ij\}\in E$,
then we have $\tau_{ij}=0$ and $K_i=\sum_{j\sim i}K_{ij}=4\pi-\sum_{\{ijkl\}\in T}\alpha_{i,jkl}$,
which is the combinatorial scalar curvature introduced by Cooper-Rivin \cite{CR}. Here $\alpha_{i,jkl}$ is the solid angle at the vertex $i$
of the tetrahedron $\{ijkl\}\in T$.
\end{remark}

We have the following rigidity for Thurston's   sphere packings with respect to the combinatorial scalar curvature in Definition \ref{def of scalar curvature}.

\begin{theorem}\label{main theorem infinitesimal rigidity for scalar curvature}
Suppose $(M, \mathcal{T})$ is a triangulated closed connected 3-manifold
with a weight $\Phi: E\rightarrow [0,\frac{\pi}{2}]$ satisfying (\ref{condition 1}) or (\ref{condition 2}).
If $\overline{r}$ is a nondegenerate Thurston's   sphere packing metric on $(M, \mathcal{T}, \Phi)$
with $K_{ij}(\overline{r})\sin^2 \Phi_{ij}\geq 0, \forall \{ij\}\in E$,
then there exists a neighborhood $U$ of $\overline{r}$ such that
the sphere packing metrics in $U$ are determined by combinatorial scalar curvatures up to scaling.
\end{theorem}

Theorem \ref{main theorem infinitesimal rigidity for scalar curvature} is referred as the local rigidity for Thurston's   sphere
packing, which means that Thurston's   sphere packing metric is locally parameterized by the combinatorial scalar curvatures.
Note that there is a subtle difference between the notion of the infinitesimal rigidity in Theorem \ref{main theorem infinitesimal rigidity}
and the notion of local rigidity in Theorem \ref{main theorem infinitesimal rigidity for scalar curvature}.
For the combinatorial scalar curvature, the local rigidity of Thurston's   sphere packings implies
the infinitesimal rigidity,
which means that Thurston's   sphere packing on $(M, \mathcal{T}, \Phi)$
can not be deformed (except by scaling) while keeping the combinatorial scalar curvatures at the vertices fixed.
For combinatorial Ricci curvature defined on the edges,
Thurston's   sphere packings have infinitesimal rigidity by Theorem \ref{main theorem infinitesimal rigidity}
and do not have local rigidity due to the dimension difference.

\begin{remark}
The condition $\sin^2 \Phi_{ij}K_{ij}(\overline{r})\geq 0$ in Theorem \ref{main theorem infinitesimal rigidity for scalar curvature}
includes the case that $K_{ij}(\overline{r})<0$ with $\Phi_{ij}=0$ for the edge $\{ij\}\in E$.
Specially, if $\Phi\equiv 0$, there is no constraints on the combinatorial Ricci curvatures
and Theorem \ref{main theorem infinitesimal rigidity for scalar curvature}
is reduced to the local rigidity of the tangential sphere packings in \cite{CR,G1,G2,R}.
Furthermore, the tangential sphere packings have global rigidity \cite{X1}.
\end{remark}

\begin{corollary}\label{coro of local rigidity}
Suppose $(M, \mathcal{T})$ is a triangulated closed connected  3-manifold
with a weight $\Phi: E\rightarrow [0,\frac{\pi}{2}]$ satisfying (\ref{condition 1}) or (\ref{condition 2}).
If $\overline{r}$ is a nondegenerate sphere packing metric on $(M, \mathcal{T}, \Phi)$
with $K_{ij}(\overline{r})\geq 0$ for any $\{ij\}\in E$, then
there exists a neighborhood $U$ of $\overline{r}$ such that
the sphere packing metrics in $U$ are determined by the combinatorial scalar curvatures up to scaling.
\end{corollary}

The basic idea to prove Theorem \ref{main theorem infinitesimal rigidity}
and Theorem \ref{main theorem infinitesimal rigidity for scalar curvature}
is using a variational principle with the discrete Hilbert-Einstein functional as the action functional.
There are two main difficulties for the proof.
The first difficulty comes from the connectedness of
the admissible space of sphere packing metrics for a tetrahedron, i.e. Theorem \ref{simply connect main theorem},
which is proved via the geometric center in Section \ref{section 2}.
The main idea of the proof comes from a new proof \cite{X3} of Bowers-Stephenson conjecture
for inversive distance circle packings on surfaces \cite{BS},
which simplifies the proofs of the Bowers-Stephenson conjecture in \cite{Guo,L2,X2}.
The second difficulty comes from the negative semi-definiteness of discrete Laplacian for Thurston's sphere packing metrics, i.e. Theorem \ref{main theorem positivity of Laplacian}, which is proved in Section \ref{section 3} with the help of
Glickenstein's second variational formula for the discrete Hilbert-Einstein functional.
The negative semi-definiteness is proved by the property of diagonally dominant matrix under the condition (\ref{condition 1}).
While under the condition (\ref{condition 2}), the negative semi-definiteness is proved by continuity of eigenvalues of
the discrete Laplacian and calculations of a 3-dimensional determinant which is quite difficult even in the tangential setting.

Note that the rigidity of Thurston's sphere packings on closed 3-manifolds are not fully understood in Theorem \ref{main theorem infinitesimal rigidity} and Theorem \ref{main theorem infinitesimal rigidity for scalar curvature}.
For example, very small, but not the same, perturbations of a tangential sphere packing are not allowed in Theorem \ref{main theorem infinitesimal rigidity} and Theorem \ref{main theorem infinitesimal rigidity for scalar curvature}, while the tangential sphere packings on
closed 3-manifolds are proved to be globally rigid \cite{X1}.
It is natural to ask whether the rigidity is true for all Thurston's sphere packing metrics on closed 3-manifolds.
Motivated by Theorem \ref{main theorem infinitesimal rigidity}, Theorem \ref{main theorem infinitesimal rigidity for scalar curvature}
and the global rigidity for the tangential sphere packings in \cite{X1}, we have the following conjecture on the rigidity of
Thurston's sphere packings on closed 3-manifolds.

\begin{conjecture}\label{conj}
Suppose $(M, \mathcal{T})$ is a triangulated closed connected  3-manifold
with a weight $\Phi: E\rightarrow [0,\frac{\pi}{2}]$.
In the space of Thurston's sphere packing metrics on  $(M, \mathcal{T}, \Phi)$ with $K_{ij}\sin^2 \Phi_{ij}\geq 0, \forall \{ij\}\in E$,
the sphere packing metric is determined by the combinatorial scalar curvature up to scaling.
\end{conjecture}

Some constraints on the weight $\Phi: E\rightarrow [0,\frac{\pi}{2}]$ as that in Theorem \ref{simply connect main theorem} may be needed in Conjecture \ref{conj}
to ensure that the admissible space of Thurston's sphere packing metrics on a tetrahedron is simply connected.
It is believed that the recent work of Doehrman-Glickenstein \cite{DG} will play a key role in the proof of Conjecture \ref{conj}, especially
in the proof of the negative semi-definiteness of discrete Laplacian for Thurston's sphere packing metrics.
%
%

The paper is organized as follows.
In Section \ref{section 2}, we prove the simply connectedness of
the admissible space of Thurston's sphere packing metrics for a tetrahedron, i.e. Theorem \ref{simply connect main theorem}.
In Section \ref{section 3}, we introduce the definition of discrete Laplacian for Thurston's sphere packing metrics and
prove its negative semi-definiteness, i.e. Theorem \ref{main theorem positivity of Laplacian}.
In Section \ref{section 4}, using Glickenstein's variational formulas for general discrete conformal variations of piecewise flat 3-manifolds
and  the negative semi-definiteness of discrete Laplacian,
we prove the main results, i.e.  Theorem \ref{main theorem infinitesimal rigidity} and Theorem \ref{main theorem infinitesimal rigidity for scalar curvature}.

~\\
\textbf{Acknowledgements}
The authors thank the referee for reading the paper carefully, and providing many valuable and helpful comments, which help to improve the paper significantly.
The work was started in March 2017 and finished when the first author was visiting University of British Columbia and the second author was visiting
Rutgers University.
The second author thanks Professor Feng Luo for his invitation to Rutgers University and communications and interesting on this work
and thanks Professor Tian Yang for communications. The first author was supported by the Grant of NSF of Hunan province no. 2018JJ2073.
The second author was supported by Hubei Provincial Natural Science Foundation of China under grant no. 2017CFB681,
Fundamental Research Funds for the Central Universities and
National Natural Science Foundation of China under grant no. 61772379 and no. 11301402.

\section{Admissible space of Thurston's sphere packing metrics for a tetrahedron}\label{section 2}
In this section, we study the admissible space of Thurston's   sphere packing metrics
for a tetrahedron $\sigma=\{1234\}\in T$. Denote the vertices set, edge set and face set of the tetrahedron $\sigma$ as $V_\sigma, E_\sigma$ and $F_\sigma$ respectively.
Set
\begin{equation*}
\begin{aligned}
G_0(l)=\left[
                   \begin{array}{ccccc}
                     0 & 1 & 1 & 1 & 1 \\
                     1 & 0 & l_{12}^2 & l_{13}^2 & l_{14}^2 \\
                     1 & l_{12}^2 & 0 & l_{23}^2 & l_{24}^2 \\
                     1 & l_{13}^2 & l_{23}^2 & 0 & l_{34}^2 \\
                     1 & l_{14}^2 & l_{24}^2 & l_{34}^2 & 0 \\
                   \end{array}
                 \right],
\end{aligned}
\end{equation*}
which is the Cayley-Menger matrix.

\begin{theorem}[\cite{Blu}]\label{nondegenerate condition}
A Euclidean tetrahedron with edge lengths $l_{12},l_{13},l_{14},l_{23},l_{24},l_{34}$ is nondegenerate in $\mathbb{R}^3$
if and only if the lengths satisfy the triangle inequalities for each face and
$\det G_0>0,$
where $\det$ is the determinant.
\end{theorem}

\begin{figure}[!htb]
\centering
  \includegraphics[height=0.4\textwidth,width=0.42\textwidth]{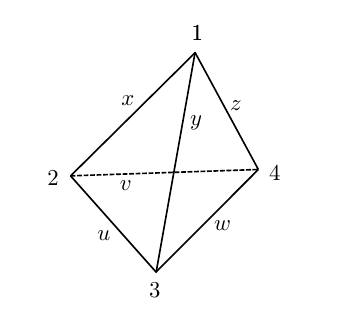}
  \caption{tetrahedron}
  \label{tetrahedron}
\end{figure}
It is proved \cite{Blu} that  $\det G_0(l)$ is a positive scalar multiplication of the square of the volume of the tetrahedron generated by $l_{12},l_{13},l_{14},l_{23},l_{24},l_{34}$.
Furthermore, the triangle inequalities on the faces are satisfied if and only if the $4\times 4$ minors of $G_0(l)$ containing the $(1,1)$-entry are negative \cite{Blu}.
To simplify the notations, we set
\begin{equation*}
\kappa_i=r_i^{-1}, \cos\Phi_{12}=x, \cos\Phi_{13}=y, \cos\Phi_{14}=z, \cos\Phi_{23}=u, \cos\Phi_{24}=v, \cos\Phi_{34}=w.
\end{equation*}
See Figure \ref{tetrahedron} for the notations. By direct calculations, we have
\begin{equation*}
\det G_0=8(r_1r_2r_3r_4)^2Q,
\end{equation*}
where
\begin{equation}\label{Q}
\begin{aligned}
Q=&\kappa_1^2(1-u^2-v^2-w^2-2uvw)+\kappa_2^2(1-y^2-z^2-w^2-2yzw)\\
&+\kappa_3^2(1-x^2-z^2-v^2-2xzv)+\kappa_4^2(1-x^2-y^2-u^2-2xyu)\\
&+2\kappa_1\kappa_2[x(1-w^2)+yu+zv+yvw+zuw]\\
&+2\kappa_1\kappa_3[y(1-v^2)+xu+zw+xvw+zuv]\\
&+2\kappa_1\kappa_4[z(1-u^2)+xv+yw+xuw+yuv]\\
&+2\kappa_2\kappa_3[u(1-z^2)+xy+vw+xzw+yzv]\\
&+2\kappa_2\kappa_4[v(1-y^2)+xz+uw+xyw+yzu]\\
&+2\kappa_3\kappa_4[w(1-x^2)+yz+uv +xyv+xzu].
\end{aligned}
\end{equation}
If the edge lengths are given by (\ref{Thurston Euclidean length introduction}) with weights in $[0, \frac{\pi}{2}]$,
the triangle inequalities on the faces are satisfied (\cite{T1} Lemma 13.7.2, \cite{X2, Z}).
By Theorem \ref{nondegenerate condition}, we have the following criteria for non-degeneracy of Thurston's sphere packing metrics.

\begin{corollary}\label{nondegenerate condition for tetrahedron}
A tetrahedron $\sigma=\{1234\}$ generated by
a Thurston's   sphere packing metric is nondegenerate in $\mathbb{R}^3$ if and only if $Q>0$.
\end{corollary}
For the following applications, set
\begin{equation}\label{A1, A2, A3, A4}
\begin{aligned}
A_1=&u^2+v^2+w^2+2uvw-1,
A_2=y^2+z^2+w^2+2yzw-1,\\
A_3=&x^2+z^2+v^2+2xzv-1,
\ A_4=x^2+y^2+u^2+2xyu-1.
\end{aligned}
\end{equation}

\begin{figure}[!htb]
\centering
  \includegraphics[height=0.4\textwidth,width=0.45\textwidth]{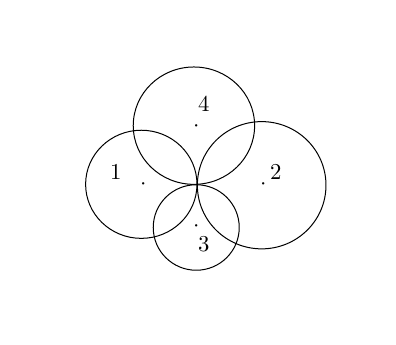}
  \caption{Degenerate sphere packing metric}
  \label{four_circle}
\end{figure}

Note that there are some trivial cases that, for any $r=(r_1, r_2,r_3, r_4)\in \mathbb{R}^4_{>0}$, the tetrahedron generated by
the sphere packing metric $r$ degenerates.
For example, if the weights on the edges of the tetrahedron satisfy
$x=w=1$ and  $y=z=u=v=0$,
then it is straight forward to check $Q=0$ for any $(r_1, r_2,r_3, r_4)\in \mathbb{R}^4_{>0}$,
which implies the tetrahedron degenerates by Corollary \ref{nondegenerate condition for tetrahedron}. Geometrically, the conditions $x=w=1$ and  $y=z=u=v=0$ mean that two pairs of tangential attached spheres meet orthogonally.
Please refer to Figure \ref{four_circle} for a configuration of the four spheres.
We can also choose another pair of opposite edges with weight $1$ and the other edges with weight $0$, in
which case the tetrahedron also degenerates for any $(r_1, r_2,r_3, r_4)\in \mathbb{R}^4_{>0}$.
We will prove that these are the only exceptional cases that
the admissible space of Thurston's sphere packing metrics for a tetrahedron is empty.
Set $S$ as the set of weights on the edges of the tetrahedron with value $1$ on a pair of opposite edges and value $0$ on the other edges, i.e.
\begin{eqnarray}\label{exception weights}
S=\{(1,0,0,0,0,1), (0,1,0,0,1,0), (0,0,1,1,0,0)\}.
\end{eqnarray}
We require the weight on the edges to be a map $\eta: E_{\sigma}\rightarrow [0,1]$ with $\eta(E_\sigma)\in W=[0, 1]^6\setminus S$ in the following, where $\eta_{ij}=\cos\Phi_{ij}$.
%
Set
\begin{equation}\label{h_i}
\begin{aligned}
h_1=&\kappa_1(1-u^2-v^2-w^2-2uvw)+\kappa_2[x(1-w^2)+yu+zv+yvw+zuw]\\
    &+\kappa_3[y(1-v^2)+xu+zw+xvw+zuv]+\kappa_4[z(1-u^2)+xv+yw+xuw+yuv],\\
h_2=&\kappa_2(1-y^2-z^2-w^2-2yzw)+\kappa_1[x(1-w^2)+yu+zv+yvw+zuw]\\
    &+\kappa_3[u(1-z^2)+xy+vw+xzw+yzv]+\kappa_4[v(1-y^2)+xz+uw+xyw+yzu],\\
h_3=&\kappa_3(1-x^2-z^2-v^2-2xzv)+\kappa_1[y(1-v^2)+xu+zw+xvw+zuv]\\
    &+\kappa_2[u(1-z^2)+xy+vw+xzw+yzv]+\kappa_4[w(1-x^2)+yz+uv +xyv+xzu],\\
h_4=&\kappa_4(1-x^2-y^2-u^2-2xyu)+\kappa_1[z(1-u^2)+xv+yw+xuw+yuv]\\
    &+\kappa_2[v(1-y^2)+xz+uw+xyw+yzu]+\kappa_3[w(1-x^2)+yz+uv +xyv+xzu].
\end{aligned}
\end{equation}

By Corollary \ref{nondegenerate condition for tetrahedron},
$r=(r_1, r_2, r_3, r_4)\in \mathbb{R}^4_{>0}$ is a degenerate Thurston's sphere packing metric if and only if the corresponding $\kappa=(\kappa_1, \kappa_2, \kappa_3, \kappa_4)\in \mathbb{R}^4_{>0}$ satisfies
\begin{equation}\label{decomposition of Q}
\begin{aligned}
Q=\kappa_1h_1+\kappa_2h_2+\kappa_3h_3+\kappa_4h_4\leq 0.
\end{aligned}
\end{equation}

\begin{remark}\label{geometric meaning of h_i}
Similar decomposition of $Q$ in (\ref{decomposition of Q}) appears in \cite{G1,G2} for tangential sphere packings.
$h_l$ in (\ref{h_i}) is closely related to the signed distance $h_{ijk,l}$ of the geometric center $c_{ijkl}$ of the tetrahedron to
the face $\{ijk\}$. More precisely, if the tetrahedron is nondegenerate, we have
\begin{equation}\label{relation of h_{ijkl} and h_4}
\begin{aligned}
h_{ijk,l}
=\frac{r_i^2r_j^2r_k^2r_l^2 }{12A_{ijk}V_{ijkl}}\kappa_lh_l,
\end{aligned}
\end{equation}
where $A_{ijk}$ is the area of the triangle $\{ijk\}$ and $V_{ijkl}$ is the volume of the tetrahedron.
Please refer to Subsection \ref{subsection definition of discrete Laplacian} for the definition of geometric center and signed distance
and refer to Subsubsection \ref{subsection simple case of negative Lap}
for the derivation of the formula (\ref{relation of h_{ijkl} and h_4}).
Note that $h_{ijk,l}$ is defined only for nondegenerate sphere packing metrics, while
$h_l$ is defined for any $r=(r_1, r_2, r_3, r_4)\in \mathbb{R}^4_{>0}$.
\end{remark}

Suppose $(r_1, r_2, r_3, r_4)\in \mathbb{R}^4_{>0}$ is a degenerate Thurston's sphere packing metric for a tetrahedron $\sigma=\{1234\}$,
which implies $Q\leq 0$ by Corollary \ref{nondegenerate condition for tetrahedron}.
In this case, we will prove none of $h_1, h_2, h_3, h_4$ is zero. Furthermore,
only one of $h_1, h_2, h_3, h_4$ is negative and the others are positive.

Note that $Q\leq 0$ could be written as a quadratic inequality of $\kappa_1$
\begin{eqnarray*}
A_1\kappa_1^2+B_1\kappa_1+C_1\geq 0
\end{eqnarray*}
with
\begin{equation}\label{A_1, B_1, C_1}
\begin{aligned}
A_1=&u^2+v^2+w^2+2uvw-1,\\
B_1=&-[2\kappa_4(vx+uwx+uvy+wy+z-zu^2)+2\kappa_3(ux+vwx+uvz+wz+y-yv^2)\\
&+2\kappa_2(x-xw^2+uy+vwy+vz+uwz)],\\
C_1=&-[\kappa_4^2(1-u^2-x^2-y^2-2uxy)+\kappa_3^2(1-v^2-x^2-z^2-2vxz)\\
&+\kappa_2^2(1-w^2-y^2-z^2-2wyz)+2\kappa_3\kappa_4(uv+w-wx^2+vxy+uxz+yz)\\
&+2\kappa_2\kappa_4(uw+wxy+v-vy^2+xz+uyz)\\
&+2\kappa_2\kappa_3(vw+xy+wxz+vyz+u-uz^2)].
\end{aligned}
\end{equation}
Lengthy direct calculations show that the discriminant $\Delta_1=B_1^2-4A_1C_1$ is
\begin{equation}
\begin{aligned}
\Delta_1
=-4\bigg{[}
 &\kappa_2^2(1-w^2)+\kappa_3^2(1-v^2)+\kappa_4^2(1-u^2)\\
 &+2\kappa_2\kappa_3(u+vw)+2\kappa_2\kappa_4(v+uw)+2\kappa_3\kappa_4(w+uv)\bigg{]}\cdot \\
 &\bigg{[}1-w^2
 -x^2-y^2-v^2-u^2-z^2+w^2x^2+v^2y^2+u^2z^2\\
&-2wyz-2vxz-2uxy-2uvw-2vxwy-2uxwz-2uvyz\bigg{]}.
\end{aligned}
\end{equation}
Note that
$$\kappa_2^2(1-w^2)+\kappa_3^2(1-v^2)+\kappa_4^2(1-u^2)
 +2\kappa_2\kappa_3(u+vw)+2\kappa_2\kappa_4(v+uw)+2\kappa_3\kappa_4(w+uv)>0$$
for any $(\kappa_2, \kappa_3, \kappa_4)\in \mathbb{R}^3_{>0}$ and $u,v,w\in [0,1]$,
which is in fact the square of the area of the triangle $\{234\}$ up to a positive factor $\frac{1}{4}r_2^2r_3^2r_4^2$ \cite{Guo,X2,Z}.
We have the following lemma on the discriminant $\Delta_1$.
\begin{lemma}\label{positive Delta}
If $A_1=u^2+v^2+w^2+2uvw-1>0$, then
the discriminant $\Delta_1>0$.
\end{lemma}
\proof
To check the positivity of the discriminant $\Delta_1$, we just need to check the negativity of the following function
\begin{equation*}
\begin{aligned}
G=&1-x^2-y^2-z^2-u^2-v^2-w^2+x^2 w^2+y^2 v^2+z^2 u^2\\
&-2wyz-2vxz-2uxy-2uvw-2vxwy-2uxwz-2uvyz.
\end{aligned}
\end{equation*}
By direct calculations, we have
\begin{equation*}
\begin{aligned}
\frac{\partial G}{\partial x}
=-2x(1-w^2)-2(yu+zv+yvw+zuw)\leq 0,
\end{aligned}
\end{equation*}
which implies
\begin{equation}\label{G|_{x=0}}
\begin{aligned}
G\leq G|_{x=0}
=&1-y^2-z^2-u^2-v^2-w^2+y^2v^2+z^2u^2-2yzw-2uvw-2yzuv.
\end{aligned}
\end{equation}
Denote the last function in (\ref{G|_{x=0}}) as $G_1$.
Taking the derivative of $G_1$ with respect to $y$ gives
\begin{eqnarray*}
\frac{\partial G_1}{\partial y}=-2y(1-v^2)-2z(w+uv)\leq 0,
\end{eqnarray*}
which implies
\begin{eqnarray*}
G_1\leq G_1|_{y=0}=-(u^2+v^2+w^2+2uvw-1)-z^2(1-u^2).
\end{eqnarray*}
By the condition $u^2+v^2+w^2+2uvw-1>0$, we have $G_1< 0$, which implies $G<0$. Therefore, the discriminant $\Delta_1>0$.\qed

\begin{remark}
If we take $Q$ as a quadratic function of $\kappa_i$, $i\in \{2,3,4\}$, results similar to Lemma \ref{positive Delta} hold for $i\in \{2,3,4\}$.
\end{remark}

\begin{lemma}\label{lemma on h not identical 0}
Suppose $(r_1, r_2, r_3, r_4)\in \mathbb{R}^4_{>0}$ is a
degenerate Thurston's sphere packing metric on the tetrahedron $\sigma=\{1234\}$
with the weight $\eta(E_\sigma)\in W=[0,1]^6\setminus S$.
Then none of $h_1, h_2, h_3, h_4$ is zero.
\end{lemma}
\proof
We prove the lemma by contradiction.
If one of $h_1, h_2, h_3, h_4$ is zero, there should be another of $h_1, h_2, h_3, h_4$ which is nonpositive by (\ref{decomposition of Q}).
Without loss of generality, we assume $h_1=0$ and $h_2\leq 0$.
By $h_1=0$, we have
\begin{equation}\label{h_1=0 in proof}
\begin{aligned}
\kappa_1(u^2+v^2+w^2+2uvw-1)=&\kappa_2[x(1-w^2)+yu+zv+yvw+zuw]\\
    &+\kappa_3[y(1-v^2)+xu+zw+xvw+zuv]\\
    &+\kappa_4[z(1-u^2)+xv+yw+xuw+yuv]\geq 0,
\end{aligned}
\end{equation}
where we have used the fact that $u,v,w\in[0,1]$ and $x,y,z$ are all nonnegative. Eq.(\ref{h_1=0 in proof})
 implies that
$A_1=u^2+v^2+w^2+2uvw-1\geq 0.$
If $u^2+v^2+w^2+2uvw-1>0$, we take $Q=\kappa_1h_1+\kappa_2h_2+\kappa_3h_3+\kappa_4h_4\leq 0$
as a quadratic inequality $A_1\kappa_1^2+B_1\kappa_1+C_1\geq 0$ of $\kappa_1$.
Then we have $\Delta_1>0$ by Lemma \ref{positive Delta}, which implies
\begin{eqnarray*}
\kappa_1\geq\frac{-B_1+\sqrt{\Delta_1}}{2A_1} \ \ \ \mbox{or} \ \ \  \kappa_1\leq\frac{-B_1-\sqrt{{\Delta_1}}}{2{A_1}}.
\end{eqnarray*}
Therefore, by the definition of $h_1$ in (\ref{h_i}), we have
\begin{eqnarray*}
h_1=-\frac{1}{2}(2A_1\kappa_1+B_1)\leq -\frac{1}{2}\sqrt{\Delta_1}<0 \ \ \ \mbox{or} \ \ \  h_1=-\frac{1}{2}(2A_1\kappa_1+B_1)\geq \frac{1}{2}\sqrt{\Delta_1}>0
\end{eqnarray*}
respectively,
which contradicts to $h_1=0$. Therefore, we have
\begin{equation}\label{equation A_1}
\begin{aligned}
u^2+v^2+w^2+2uvw-1=0,
\end{aligned}
\end{equation}
which implies
\begin{equation}\label{equation c_ij}
\begin{aligned}
&x(1-w^2)+y(u+vw)+z(v+uw)=0,\\
&y(1-v^2)+x(u+vw)+z(w+uv)=0,\\
&z(1-u^2)+x(v+uw)+y(w+uv)=0
\end{aligned}
\end{equation}
by (\ref{h_1=0 in proof}). By (\ref{equation A_1}), at least one of $u, v, w$ is positive.
Without loss of generality, we assume $u>0$.
Then we have $x=y=0$ and $z(1-u^2)=0$ by (\ref{equation c_ij}), which implies $z=0$ or $u=1$.

Note that $h_2\leq 0$  is equivalent to
\begin{equation}\label{equation h_2<0}
\begin{aligned}
\kappa_2(y^2+z^2+w^2+2yzw-1)\geq &\kappa_1[x(1-w^2)+u(y+zw)+v(z+yw)]\\
    &+\kappa_3[u(1-z^2)+v(w+yz)+x(y+wz)]\\
    &+\kappa_4[v(1-y^2)+x(z+yw)+u(w+yz)].
\end{aligned}
\end{equation}
In the case $x=y=z=0$, (\ref{equation h_2<0}) is equivalent to
\begin{equation*}
\begin{aligned}
0\geq (w^2-1)\kappa_2\geq\kappa_3(u+vw)+\kappa_4(v+uw)>0,
\end{aligned}
\end{equation*}
which is impossible.
In the case $x=y=0, u=1$, we have $v=w=0$ by (\ref{equation A_1}). Then (\ref{equation h_2<0}) is equivalent to
\begin{equation*}
\begin{aligned}
0\geq (z^2-1)\kappa_2\geq\kappa_3(1-z^2)\geq 0,
\end{aligned}
\end{equation*}
which implies $z=1$. This contradicts the condition $\eta(E_\sigma)\in W=[0,1]^6\setminus S$.
This completes the proof.\qed

\begin{remark}
Lemma \ref{lemma on h not identical 0} implies that, for a degenerate Thurston's sphere packing metric on a tetrahedron $\sigma$ with the weight $\eta(E_\sigma)\in W=[0,1]^6\setminus S$, the geometric center
can never lie in the planes determined by the faces. 
\end{remark}

By (\ref{decomposition of Q}), Corollary \ref{nondegenerate condition for tetrahedron} and Lemma \ref{lemma on h not identical 0}, at least one of $h_1, h_2, h_3, h_4$ is negative and
the others are nonzero if $(r_1, r_2, r_3, r_4)\in \mathbb{R}^4_{>0}$ is a
degenerate Thurston's sphere packing metric on the tetrahedron.
Furthermore, we have the following result.

\begin{lemma}\label{lemma on h1<0 h2<0 at the same time}
Suppose $(r_1, r_2, r_3, r_4)\in \mathbb{R}^4_{>0}$ is
a degenerate Thurston's   sphere packing metric on a tetrahedron with $\eta(E_\sigma)\in W=[0,1]^6\setminus S$.
Then there exists no subset $\{i,j\}\subset \{1,2,3,4\}$ such that $h_i<0$ and $h_j<0$.
\end{lemma}
\proof
Without loss of generality, we assume $h_1< 0$ and $h_2< 0$.
Then we have
\begin{equation*}
\begin{aligned}
&\kappa_1(u^2+v^2+w^2+2uvw-1)>\kappa_2[x(1-w^2)+yu+zv+yvw+zuw]\\
    &+\kappa_3[y(1-v^2)+xu+zw+xvw+zuv]+\kappa_4[z(1-u^2)+xv+yw+xuw+yuv]
\end{aligned}
\end{equation*}
and
\begin{equation*}
\begin{aligned}
&\kappa_2(y^2+z^2+w^2+2yzw-1)>\kappa_1[x(1-w^2)+yu+zv+yvw+zuw]\\
    &+\kappa_3[u(1-z^2)+xy+vw+xzw+yzv]+\kappa_4[v(1-y^2)+xz+uw+xyw+yzu],
\end{aligned}
\end{equation*}
which implies
\begin{equation*}
\begin{aligned}
u^2+v^2+w^2+2uvw-1>0, y^2+z^2+w^2+2yzw-1>0,
\end{aligned}
\end{equation*}
and the following function
\begin{equation*}
\begin{aligned}
F:=&(u^2+v^2+w^2+2uvw-1)(y^2+z^2+w^2+2yzw-1)\\
&-[x(1-w^2)+yu+zv+yvw+zuw]^2
\end{aligned}
\end{equation*}
is positive.

We claim that $F\leq 0$ under the condition $x,y,z,u,v,w\in [0,1]$ and $u^2+v^2+w^2+2uvw-1>0$.
Then the proof of the lemma follows by contradiction.
The proof of the claim is as follows.

$F$ is decreasing in $x$, which is obvious. Therefore, we have
\begin{equation}\label{F|_{x=0}}
\begin{aligned}
F\leq F|_{x=0}
=&(u^2+v^2+w^2+2uvw-1)(y^2+z^2+w^2+2yzw-1)\\
 &-[y(u+vw)+z(v+uw)]^2.
\end{aligned}
\end{equation}
Denote the last function in (\ref{F|_{x=0}}) of $y,z,u,v,w$ as $F_1$.
By direct calculations, we have
\begin{equation*}
\begin{aligned}
\frac{1}{2}\frac{\partial F_1}{\partial y}
=&(w^2-1)[y(1-v^2)+z(w+uv)]\leq 0.
\end{aligned}
\end{equation*}
Therefore
\begin{equation}\label{F_1|_{y=0}}
\begin{aligned}
F_1\leq F_1|_{y=0}=(u^2+v^2+w^2+2uvw-1)(z^2+w^2-1)-[z(v+uw)]^2.
\end{aligned}
\end{equation}
Denote the last function in (\ref{F_1|_{y=0}}) as $F_2$.
Note that
$
\frac{\partial F_2}{\partial z}
=-2z(1-u^2)(1-w^2)\leq 0.
$
Therefore
\begin{equation*}
\begin{aligned}
F_2\leq F_2|_{z=0}=(u^2+v^2+w^2+2uvw-1)(w^2-1)\leq 0.
\end{aligned}
\end{equation*}
In summary, we have
\begin{equation*}
\begin{aligned}
F\leq F|_{x=0}=F_1\leq F_1|_{y=0}=F_2\leq F_2|_{z=0}\leq 0.
\end{aligned}
\end{equation*}
This completes the proof of the claim. \qed

Combining Lemma  \ref{lemma on h not identical 0} and Lemma \ref{lemma on h1<0 h2<0 at the same time}, we have the following proposition.

\begin{proposition}\label{proposition h_i<0 others positive}
Suppose $(r_1, r_2, r_3, r_4)\in \mathbb{R}^4_{>0}$ is
a degenerate Thurston's   sphere packing metric on a tetrahedron with $\eta(E_\sigma)\in W=[0,1]^6\setminus S$.
Then one of $h_1, h_2, h_3, h_4$ is negative and the others are positive.
\end{proposition}

\begin{remark}
Combining with Remark \ref{geometric meaning of h_i},
Proposition \ref{proposition h_i<0 others positive} has an interesting geometric explanation as follows.
Suppose a nondegenerate tetrahedron $\sigma=\{1234\}$ becomes degenerate as Thurston's sphere packing metric $r=(r_1, r_2, r_3, r_4)\rightarrow \bar{r}=(\bar{r}_1, \bar{r}_2, \bar{r}_3, \bar{r}_4)\in \partial (\Omega_{1234}(\eta))$,
where $\Omega_{1234}(\eta)\subset \mathbb{R}^4_{>0}$ is
the admissible space of Thurston's
sphere packing metrics for a given weight $\eta$ with $\eta(E_\sigma)\in W=[0,1]^6\setminus S$ and
$\partial (\Omega_{1234}(\eta))$ is the boundary of $\Omega_{1234}(\eta)$ in $\mathbb{R}^4_{>0}$.
Then the volume $V_{1234}$ of the tetrahedron $\{1234\}$ goes to zero as $r\rightarrow\bar{r}\in \partial (\Omega_{1234}(\eta))$.
Without loss of generality, we can assume $h_1(\bar{r})<0, h_2(\bar{r})>0,h_3(\bar{r})>0,h_4(\bar{r})>0$ by Proposition \ref{proposition h_i<0 others positive}.
By Remark \ref{geometric meaning of h_i}, we have $h_{234,1}\rightarrow -\infty, h_{123,4}\rightarrow +\infty, h_{124,3}\rightarrow +\infty, h_{134,2}\rightarrow +\infty$.
This implies that, as $r\rightarrow\bar{r}\in \partial (\Omega_{1234}(\eta))$,  the geometric center $c_{1234}$ and the tetrahedron $\{1234\}$  will lie in different half spaces relative to the plane determined by the face $\{234\}$, while $c_{1234}$ and the tetrahedron $\{1234\}$  will lie in the same half spaces relative to the planes determined by the faces $\{123\}, \{124\}, \{134\}$.
\end{remark}

Now we can prove the main result of this section.
\begin{theorem}\label{simply connect main theorem}
Suppose $\sigma=\{1234\}$ is a tetrahedron of $(M, \mathcal{T})$. Then the admissible space $\Omega_{1234}(\eta)$ of Thurston's
sphere packing metrics for a given weight $\eta$ with $\eta(E_\sigma)\in W=[0,1]^6\setminus S$
is a simply connected nonempty set.
\end{theorem}
\proof
Suppose $r=(r_1, r_2, r_3, r_4)$ is a degenerate sphere packing metric, then the corresponding
$\kappa=(\kappa_1, \kappa_2, \kappa_3, \kappa_4)\in \mathbb{R}^4_{>0}$
satisfies (\ref{decomposition of Q}).
By Proposition \ref{proposition h_i<0 others positive},
one of $h_1, h_2, h_3, h_4$ is negative and the others are positive.
Without loss of generality, we assume $h_1<0$.
Then we have
$A_1=u^2+v^2+w^2+2uvw-1>0$ by the definition of $h_1$.
Take $Q\leq 0$ as a quadratic inequality $A_1\kappa_1^2+B_1\kappa_1+C_1\geq 0$ of $\kappa_1$,
where $A_1, B_1, C_1$ are given by (\ref{A_1, B_1, C_1}).
As $A_1=u^2+v^2+w^2+2uvw-1>0$, we have $\Delta_1=B_1^2-4A_1C_1>0$ by Lemma \ref{positive Delta}.
Therefore,
\begin{eqnarray*}
\kappa_1\geq\frac{-B_1+\sqrt{\Delta_1}}{2A_1}\  \mbox{or} \ \
\kappa_1\leq\frac{-B_1-\sqrt{\Delta_1}}{2A_1}.
\end{eqnarray*}
Note that $h_1< 0$ is equivalent to $\kappa_1> \frac{-B_1}{2A_1}$ by definition of $h_1$.
Therefore,
$\kappa_1\geq\frac{-B_1+\sqrt{\Delta_1}}{2A_1}$,
which implies $\kappa_1>\frac{-B_1}{2A_1}$ and $h_1< 0$.

In the case $A_1=u^2+v^2+w^2+2uvw-1>0$, we set
\begin{eqnarray*}
U_1=\{((r_1, r_2, r_3, r_4)\in \mathbb{R}^4_{>0}|\kappa_1\geq \frac{-B_1+\sqrt{\Delta_1}}{2A_1}\},
\end{eqnarray*}
which is a closed domain in $\mathbb{R}^4_{>0}$ bounded by an analytic graph on $\mathbb{R}^3_{>0}$ and contained in $\mathbb{R}^4_{>0}\setminus \Omega_{1234}(\eta)$.
$U_2, U_3, U_4$ are defined similarly, if the corresponding $A_i>0$.
Note that $U_1, U_2, U_3, U_4$ are mutually disjoint by Lemma \ref{lemma on h1<0 h2<0 at the same time}, if they are nonempty.

If $A_\mu>0$ for $\mu\in P\subset \{1,2,3,4\}$ and $A_\nu\leq 0$ for $\nu\in \{1,2,3,4\}\setminus P$, then the space of degenerate sphere packing metrics
is $\cup_{\mu\in P}U_\mu$.
The corresponding admissible space of sphere packing metric is
\begin{equation*}
\begin{aligned}
\Omega_{1234}(\eta)=\mathbb{R}^4_{>0}\setminus\cup_{\mu\in P}U_\mu,
\end{aligned}
\end{equation*}
which is nonempty.
Note that $\Omega_{1234}(\eta)$ is homotopy equivalent to $\mathbb{R}^4_{>0}$ because
$U_1, U_2, U_3, U_4$ are mutually disjoint and bounded by analytic graphs on $\mathbb{R}^3_{>0}$.
This implies that, for any fixed weight $\eta$ with $\eta(E_\sigma)\in W=[0, 1]^6\setminus S$,
the admissible space $\Omega_{1234}(\eta)$ of sphere packing metrics is simply connected.  \qed

\begin{corollary}\label{coro of adm space is R4}
Suppose $\sigma=\{1234\}$ is a tetrahedron of $(M, \mathcal{T})$.
$\eta$ is a weight on $E_\sigma$ with $\eta(E_\sigma)\in W=[0,1]^6\setminus S$.
Then the admissible space $\Omega_{1234}(\eta)=\mathbb{R}^4_{>0}$ if and only if all of the inequalities
\begin{equation*}
\begin{aligned}
A_1=&u^2+v^2+w^2+2uvw-1\leq 0,
A_2=y^2+z^2+w^2+2yzw-1\leq 0,\\
A_3=&x^2+z^2+v^2+2xzv-1\leq 0,
\ A_4=x^2+y^2+u^2+2xyu-1\leq 0.
\end{aligned}
\end{equation*}
are true.
Specially, if
\begin{equation}\label{sum of Phi > pi}
\Phi_{ij}+\Phi_{ik}+\Phi_{jk}\geq \pi
\end{equation}
for any face $\{ijk\}\in F_\sigma$, then the admissible space $\Omega_{1234}(\eta)$ of Thurston's   sphere packing metrics is $\mathbb{R}^{4}_{>0}$.
\end{corollary}

\proof
The first part of Corollary \ref{coro of adm space is R4} follows from the proof of Theorem \ref{simply connect main theorem}.
For the second part, note that
\begin{equation*}
\begin{aligned}
A_1
=&\cos^2\Phi_{23}+\cos^2\Phi_{24}+\cos^2\Phi_{34}+2\cos\Phi_{23}\cos\Phi_{24}\cos\Phi_{34}-1\\
=&(\cos\Phi_{23}+\cos\Phi_{24}\cos\Phi_{34})^2-\sin^2\Phi_{24}\sin^2\Phi_{34}\\
=&4\cos\frac{\Phi_{23}+\Phi_{24}+\Phi_{34}}{2}\cos\frac{\Phi_{24}+\Phi_{34}-\Phi_{23}}{2}
\cos\frac{\Phi_{23}-\Phi_{24}+\Phi_{34}}{2}\cos\frac{\Phi_{23}+\Phi_{24}-\Phi_{34}}{2}.
\end{aligned}
\end{equation*}
For $\Phi_{ij}\in [0, \frac{\pi}{2}]$, we have
\begin{equation*}
\begin{aligned}
\cos\frac{\Phi_{24}+\Phi_{34}-\Phi_{23}}{2}\geq 0,
\cos\frac{\Phi_{23}-\Phi_{24}+\Phi_{34}}{2}\geq 0,
\cos\frac{\Phi_{23}+\Phi_{24}-\Phi_{34}}{2}\geq 0.
\end{aligned}
\end{equation*}
Therefore, if $\Phi_{23}+\Phi_{24}+\Phi_{34}\geq  \pi,$ we have
$\cos\frac{\Phi_{23}+\Phi_{24}+\Phi_{34}}{2}\leq 0$ and $A_1\leq 0$.
Similarly, if
$$\Phi_{13}+\Phi_{14}+\Phi_{34}\geq \pi, \Phi_{12}+\Phi_{14}+\Phi_{24}\geq \pi, \Phi_{12}+\Phi_{13}+\Phi_{23}\geq \pi,$$
we have $A_2\leq 0, A_3\leq0, A_4\leq0$ respectively. \qed

\begin{remark}
The condition (\ref{sum of Phi > pi}) in Corollary \ref{coro of adm space is R4}
includes the case that $\Phi_{ij}\in [\frac{\pi}{3}, \frac{\pi}{2}]$.
Furthermore, the condition (\ref{sum of Phi > pi}) allows some of the intersection angles to take
values in $(0, \frac{\pi}{3})$.
However, the cases illustrated in Figure \ref{four_circle} should be excluded.
\end{remark}

\begin{remark}
The simply connectedness of the admissible space of tangential sphere packing metrics
for a tetrahedron was first proved by Cooper-Rivin \cite{CR}.
The second author \cite{X1} gave a different proof of the simply connectedness with an explicit description of
the boundary of the admissible space.
Ge-Jiang-Shen \cite{GJS} and Ge-Hua \cite{GH}
gave some new viewpoints on the method in \cite{X1} for tangential sphere packings.
The proof of Theorem \ref{simply connect main theorem} for Thurston's sphere packing metrics
unifies the proofs for simply connectedness in \cite{GJS,X1}.
This method was applied to prove the simply connectedness of the admissible space of Thurston's hyperbolic sphere packing metrics for
a tetrahedron \cite{HX1}.
This method has also been applied to give a new proof of the Bowers-Stephenson conjecture for the inversive distance circle packings
on triangulated surfaces \cite{X3}, which simplifies the proofs of the Bowers-Stephenson conjecture in \cite{Guo,L2,X2}.
\end{remark}

Using $\Omega_{1234}(\eta)$, we can further define the following 10-dimensional set
\begin{equation*}
\begin{aligned}
\Omega_{1234}=\mathop{\cup}\limits_{\eta\in W}(\eta\times\Omega_{1234}(\eta))=\{(\eta, r)\in W\times \mathbb{R}^4_{>0}|Q(\eta, r)>0\},
\end{aligned}
\end{equation*}
where we take $Q$ as a function of $\eta$ and $r$.

\begin{corollary}\label{connectness of 10-dim admissible space}
The set $\Omega_{1234}$ is connected.
\end{corollary}
The proof of Corollary \ref{connectness of 10-dim admissible space} follows easily from
Theorem \ref{simply connect main theorem}, the connectedness of $W$ together with the continuity of $Q$
and is almost the same as that of Lemma 2.8 in \cite{X3}, so we omit the proof of Corollary \ref{connectness of 10-dim admissible space} here.

If $\widetilde{W}\subset W$ is a connected subset of $W$, we can also define
\begin{equation}\label{Omega-tilde}
\begin{aligned}
\widetilde{\Omega}_{1234}=
\mathop{\cup}\limits_{\eta\in \widetilde{W}}(\eta\times\Omega_{1234}(\eta))=\{(\eta, r)\in \widetilde{W}\times \mathbb{R}^4_{>0}|Q(\eta, r)>0\}.
\end{aligned}
\end{equation}
Similar to Corollary \ref{connectness of 10-dim admissible space}, we have

\begin{corollary}\label{general connectness of 10-dim admissible space}
The set $\widetilde{\Omega}_{1234}$ is connected.
\end{corollary}

\section{Discrete Laplacian of Thurston's sphere packing metrics}\label{section 3}
\subsection{Definition of Discrete Laplacian}\label{subsection definition of discrete Laplacian}

If a nondegenerate tetrahedron is generated by a Thurston's   sphere packing metric,
then the tetrahedron can be assigned a geometric center and a geometric dual structure as follows.
Suppose $\{ijkl\}$ is a nondegenerate Euclidean tetrahedron generated by a sphere packing metric.
Then every pair of spheres attached to the vertices $s, t\in \{i,j,k,l\}$ determines a unique plane $P_{st}$ perpendicular to the edge $\{st\}$.
$P_{st}$ is the plane determined by the circle $S_s\cap S_t$ if $S_s\cap S_t\neq \emptyset$,
otherwise it is the common tangent plane of $S_s$ and $S_t$ if the two
spheres are externally tangent.
Then the six planes associated to the six edges intersect at a common
point $c_{ijkl}$ \cite{G4}, which is called the \textbf{\emph{geometric center}} of the tetrahedron $\{ijkl\}$.
Projections of the center $c_{ijkl}$ to the four faces generate the four centers $c_{ijk}, c_{ijl}, c_{ikl}, c_{jkl}$ of the four faces
$\{ijk\}, \{ijl\}, \{ikl\},\{jkl\}$ respectively.
The centers of the faces could be further projected to the edges to generate the edge centers $c_{ij}, c_{ik}, c_{il}, c_{jk}, c_{jl}, c_{kl}$
of the edges $\{ij\}, \{ik\}, \{il\}, \{jk\}, \{jl\}, \{kl\}$ respectively,
which are the intersections of edges and the corresponding perpendicular planes.
Note that the geometric center $c_{ijkl}$ may not be in the tetrahedron $\{ijkl\}$.
The signed distance of $c_{ijkl}$ to $c_{ijk}$ is denoted by $h_{ijk,l}$, which is positive if $c_{ijkl}$ is on the same side of the plane determined
by the face $\{ijk\}$ as the tetrahedron, otherwise it is negative (or zero if the center is in the plane determined by $\{ijk\}$).
Similarly, we can define the signed distance $h_{ij,k}$ of the center $c_{ijk}$ to the edge $ij$.
The signed distance of $c_{ij}$ to vertex $i$ is $d_{ij}$ and the signed distance of $c_{ij}$ to vertex $j$ is $d_{ji}$.
Obviously, we have $d_{ij}+d_{ji}=l_{ij}$.
For Thurston's sphere packing metrics,
\begin{equation*}
d_{ij}=\frac{r_i(r_i+r_j\cos\Phi_{ij})}{l_{ij}}>0.
\end{equation*}
The dual of the edge $\{ij\}$ in the tetrahedron $\{ijkl\}$ is defined to be
the planar quadrilateral determined by $c_{ij}, c_{ijk}, c_{ijkl}, c_{ijl}$
and the dual area $A_{ij,kl}$  of the edge $\{ij\}$ in the tetrahedron $\{ijkl\}$ is defined to be
the signed area of
the planar quadrilateral determined by $c_{ij}, c_{ijk}, c_{ijkl}, c_{ijl}$, i.e.
\begin{equation}\label{Aijklnew}
A_{ij,kl}=\frac{1}{2}(h_{ij,k}h_{ijk,l}+h_{ij,l}h_{ijl,k}).
\end{equation}
For a weighted triangulated 3-manifold $(M, \mathcal{T}, \Phi)$ with a nondegenerate Thurston's   sphere packing metric,
we can also define the dual area $l_{ij}^*$ of the edge $\{ij\}\in E$, which is
\begin{equation}\label{dual of an edge}
l_{ij}^*=\sum_{k,l}A_{ij,kl}=\frac{1}{2}\sum_{k,l}(h_{ij,k}h_{ijk,l}+h_{ij,l}h_{ijl,k}),
\end{equation}
where the summation is taken over the tetrahedra with $\{ij\}$ as a common edge.
Please refer to \cite{G1,G2,G3,G4,GT} for more information on geometric center and geometric dual in general cases.
Using the geometric dual, the discrete Laplacian for Thurston's   sphere packing metrics is defined as follows.
\begin{definition}[\cite{G3}]\label{definition of Laplacian}
Suppose $(M, \mathcal{T})$ is a triangulated 3-manifold with a weight $\Phi: E\rightarrow [0,\frac{\pi}{2}]$.
The discrete Laplacian of a nondegenerate sphere packing metric
is defined to be a linear map $\triangle: \mathbb{R}^V\rightarrow \mathbb{R}^V$ such that
\begin{equation*}
\triangle f_i=\sum_{j\sim i}\frac{l_{ij}^*}{l_{ij}}(f_j-f_i)
\end{equation*}
for any $f\in \mathbb{R}^V$, where $l_{ij}^*$ is the dual area of the edge $\{ij\}$ given by (\ref{dual of an edge}).
\end{definition}
The discrete Laplacian could be defined for more general discrete conformal geometric structures \cite{G3,G4} and graphs \cite{C}.
\subsection{Negative semi-definiteness of the discrete Laplacian}

We will prove that the discrete Laplacian is
negative semi-definite for a large class of Thurston's   sphere packing metrics on 3-manifolds.

\begin{theorem}\label{main theorem positivity of Laplacian}
Suppose $(M, \mathcal{T})$ is a triangulated closed connected 3-manifold with a weight $\Phi: E\rightarrow [0,\frac{\pi}{2}]$
satisfying (\ref{condition 1}) or (\ref{condition 2}).
Then the discrete Laplacian for Thurston's   sphere packing metrics
is negative semi-definite with one dimensional kernel $\{t(1, \cdots, 1)|t\in \mathbb{R}\}$.
\end{theorem}
Under the condition (\ref{condition 1}), the negative semi-definiteness of the discrete Laplacian
is proved using the property of diagonally dominant matrix.
While under the condition (\ref{condition 2}), the negative semi-definiteness of the discrete Laplacian
is proved by the continuity of the eigenvalue of the discrete Laplacian on a connected domain.
As the proofs of negative semi-definiteness of the discrete Laplacian  under the condition (\ref{condition 1}) and condition (\ref{condition 2})
are different, we present the proofs separately.

\subsubsection{Proof of Theorem \ref{main theorem positivity of Laplacian} under the condition (\ref{condition 1})}\label{subsection simple case of negative Lap}
\proof
By the Definition \ref{definition of Laplacian} of discrete Laplacian, for $f\in \mathbb{R}^{V}$, we have
\begin{equation*}
\sum_{i\in V}f_i\triangle f_i=\sum_{i\sim j}\frac{l_{ij}^*}{l_{ij}}(f_j-f_i)f_i=-\frac{1}{2}\sum_{i\sim j}\frac{l_{ij}^*}{l_{ij}}(f_j-f_i)^2.
\end{equation*}
The negative semi-definiteness of the discrete Laplacian is equivalent to
\begin{equation*}
\sum_{i\sim j}\frac{l_{ij}^*}{l_{ij}}(f_j-f_i)^2\geq 0, \ \ \ \forall f\in \mathbb{R}^{V}.
\end{equation*}
An approach to prove the negative semi-definiteness of $\triangle$ is to prove $l_{ij}^*>0$, $\forall \{ij\}\in E$.
Recall that
\begin{equation*}
l_{ij}^*=\sum_{k,l}A_{ij,kl}=\frac{1}{2}\sum_{k,l}(h_{ij,k}h_{ijk,l}+h_{ij,l}h_{ijl,k}).
\end{equation*}
To prove $l_{ij}^*>0$, we just need to prove $h_{ij,k}>0$, $h_{ij,l}>0$, $h_{ijk,l}>0$ and $h_{ijl,k}>0$.
Note that
$h_{ij,k}$ and $h_{ijk,l}$ can be computed by
\begin{equation}\label{h_{ij,k}}
h_{ij,k}=\frac{d_{ik}-d_{ij}\cos\gamma_{ijk}}{\sin\gamma_{ijk}}
\end{equation}
and
\begin{equation}\label{h_{ijk,l}}
h_{ijk,l}=\frac{h_{ij,l}-h_{ij,k}\cos\beta_{ij, kl}}{\sin \beta_{ij,kl}},
\end{equation}
where $\gamma_{ijk}$ is the inner angle at the vertex $i$ of the triangle $\{ijk\}$
and $\beta_{ij,kl}$ is the dihedral angle along the edge $\{ij\}$ in the tetrahedron $\{ijkl\}$.
See \cite{G3,G4} for (\ref{h_{ij,k}}) and (\ref{h_{ijk,l}}).

To simplify the notations, we denote
$\eta_{ij}=\cos\Phi_{ij}.$
Submitting $d_{ij}=\frac{r_i(r_i+r_j\eta_{ij})}{l_{ij}}$ and  $d_{ik}=\frac{r_i(r_i+r_k\eta_{ik})}{l_{ik}}$ into (\ref{h_{ij,k}}),
we have
\begin{equation}\label{calculation of h_ij,k}
\begin{aligned}
h_{ij,k}
=&\frac{l_{ij}r_i(r_i+r_k\eta_{ik})-l_{ik}r_i(r_i+r_j\eta_{ij})\cos\gamma_{ijk}}{l_{ij}l_{ik}\sin\gamma_{ijk}}\\
=&\frac{r_i}{2l_{ij}\cdot l_{ij}l_{ik}\sin\gamma_{ijk}}[2l_{ij}^2(r_i+r_k\eta_{ik})-2(r_i+r_j\eta_{ij})l_{ij}l_{ik}\cos\gamma_{ijk}]\\
=&\frac{r_i}{2l_{ij}\cdot l_{ij}l_{ik}\sin\gamma_{ijk}}[2l_{ij}^2(r_i+r_k\eta_{ik})-(r_i+r_j\eta_{ij})(l_{ij}^2+l_{ik}^2-l_{jk}^2)]\\
=&\frac{r_i}{l_{ij}\cdot l_{ij}l_{ik}\sin\gamma_{ijk}}[r_ir_j^2(1-\eta_{ij}^2)+r_j^2r_k(\eta_{ik}+\eta_{ij}\eta_{jk})+r_ir_jr_k(\eta_{jk}+\eta_{ij}\eta_{ik})]\\
=&\frac{1}{2l_{ij}A_{ijk}}[r_i^2r_j^2(1-\eta_{ij}^2)+r_i^2r_jr_k\Gamma_{ijk}+r_ir_j^2r_k\Gamma_{jik}],
\end{aligned}
\end{equation}
where $A_{ijk}=\frac{1}{2}l_{ij}l_{ik}\sin\gamma_{ijk}$ is the area of the triangle $\{ijk\}$, $\Gamma_{ijk}=\eta_{jk}+\eta_{ij}\eta_{ik}$ and $\Gamma_{jik}=\eta_{ik}+\eta_{ij}\eta_{jk}$.
Note that this is the result obtained in Lemma 2.5 of \cite{X2} with a different definition of $A_{ijk}$. Here we compute it directly.

It is straight forward to check that $h_{ij,k}\geq 0$ by the condition $\eta_{ij}, \eta_{ik}, \eta_{jk}\in [0, 1]$.
Furthermore, $h_{ij,k}=0$ if and only if $\eta_{ij}=1$ and $\eta_{ik}=\eta_{jk}=0$, which
corresponds to $\Phi_{ij}=0$ and $\Phi_{ik}=\Phi_{jk}=\frac{\pi}{2}$.
Therefore, under the condition (\ref{condition 1}), we have $h_{ij,k}>0$.

To compute $h_{ijk,l}$, take a unit sphere at $i$,
which determines a spherical triangle with edge lengths $\gamma_{ijk}$, $\gamma_{ijl}$, $\gamma_{ikl}$
and opposite inner angles $\beta_{il,jk}$, $\beta_{ik,jl}$, $\beta_{ij,kl}$ respectively.
By the spherical cosine law, we have
\begin{equation}\label{cos beta_ij,kl}
\begin{aligned}
\cos\beta_{ij}
=&\frac{\cos\gamma_{ikl}-\cos\gamma_{ijk}\cos\gamma_{ijl}}{\sin\gamma_{ijk}\sin\gamma_{ijl}}\\
=&\frac{l_{ij}^2l_{ik}l_{il}}{4A_{ijk}A_{ijl}}\left(\cos\gamma_{ikl}-\cos\gamma_{ijk}\cos\gamma_{ijl}\right)\\
=&\frac{l_{ij}^2l_{ik}l_{il}}{4A_{ijk}A_{ijl}}\left(\frac{l_{ik}^2+l_{il}^2-l_{kl}^2}{2l_{ik}l_{il}}
-\frac{l_{ij}^2+l_{ik}^2-l_{jk}^2}{2l_{ij}l_{ik}}\cdot \frac{l_{ij}^2+l_{il}^2-l_{jl}^2}{2l_{ij}l_{il}}\right)\\
=&\frac{1}{16A_{ijk}A_{ijl}}\left[2l_{ij}^2(l_{ik}^2+l_{il}^2-l_{kl}^2)-(l_{ij}^2+l_{ik}^2-l_{jk}^2)(l_{ij}^2+l_{il}^2-l_{jl}^2)\right],
\end{aligned}
\end{equation}
where $A_{ijk}$ and $A_{ijl}$ are the areas of triangles $\{ijk\}$ and $\{ijl\}$ respectively.
Submitting (\ref{Thurston Euclidean length introduction}), (\ref{calculation of h_ij,k}), (\ref{cos beta_ij,kl}) into (\ref{h_{ijk,l}}),
after lengthy calculations, we have
\begin{equation}\label{calculation of h_ijk,l}
\begin{aligned}
h_{ijk,l}
=&\frac{r_i^2r_j^2r_k^2r_l^2 }{12A_{ijk}V_{ijkl}}\{r_l^{-2}\left(1-\eta_{ij}^2-\eta_{ik}^2-\eta_{jk}^2-2 \eta_{ij} \eta_{ik} \eta_{jk}\right)\\
&+r_k^{-1}r_l^{-1} \left[\eta_{kl}(1-\eta_{ij}^2)+\eta_{ik}(\eta_{il}+\eta_{ij} \eta_{jl})+\eta_{jk} (\eta_{jl}+\eta_{ij} \eta_{il})\right]\\
&+r_j^{-1}r_l^{-1} \left[\eta_{jl}(1-\eta_{ik}^2)+\eta_{ij} (\eta_{il}+\eta_{ik} \eta_{kl})+\eta_{jk} (\eta_{kl}+\eta_{ik} \eta_{il})\right]\\
&+r_i^{-1}r_l^{-1} \left[\eta_{il} (1-\eta_{jk}^2)+\eta_{ij} (\eta_{jl}+\eta_{jk} \eta_{kl})+\eta_{ik} (\eta_{kl}+\eta_{jk} \eta_{jl})\right]\},
\end{aligned}
\end{equation}
where $V_{ijkl}$ is the volume of the tetrahedron and
the formula $A_{ijk}A_{ijl}\sin \beta_{ij,kl}=\frac{3}{2}l_{ij}V_{ijkl}$ is used in the calculations.
%
Similar to the proof of Corollary \ref{coro of adm space is R4}, we have $h_{ijk,l}>0$ under the condition (\ref{condition 1}).
This completes the proof of Theorem  under the condition (\ref{condition 1}). \qed

%

\subsubsection{Proof of Theorem \ref{main theorem positivity of Laplacian} under the condition (\ref{condition 2})}
For a generic weight $\Phi: E\rightarrow [0, \frac{\pi}{2}]$,
$h_{ijk,l}$ could be negative or zero. See \cite{G2} for an example in the tangential case.
The method in Subsubsection \ref{subsection simple case of negative Lap} no longer works.
In this case, we take the discrete Laplacian as a matrix $\Lambda$, where
\begin{equation*}
\begin{aligned}
\Lambda_{ij}=\left\{
               \begin{array}{ll}
                 \frac{l_{ij}^*}{l_{ij}}, & \hbox{if $j\sim i$;} \\
                 -\sum_{k\sim i}\Lambda_{ik}, & \hbox{if $j=i$;} \\
                 0, & \hbox{otherwise.}
               \end{array}
             \right.
\end{aligned}
\end{equation*}
Note that
$\Lambda=\sum_{\sigma\in T}(\Lambda_{\sigma})_{E},$
where $\Lambda_{\sigma}$ is a $4\times 4$ matrix defined for a tetrahedron $\sigma=\{ijkl\} \in T$ with
\begin{equation*}
\begin{aligned}
(\Lambda_{\sigma})_{ij}
=
\left\{
  \begin{array}{ll}
    \frac{A_{ij,kl}}{l_{ij}}, & \hbox{if $j\sim i$,} \\
    -\sum_{j\sim i}(\Lambda_{\sigma})_{ij}, & \hbox{if $j=i$,}
  \end{array}
\right.
\end{aligned}
\end{equation*}
and $(\Lambda_{\sigma})_{E}$ is the matrix
extended by zeroes to a $|V|\times |V|$ matrix so that the matrix $(\Lambda_{\sigma})_{E}$
acts on a vector $(f_1, \cdots, f_{|V|})$ only on the coordinates corresponding to the vertices in $\sigma$.
This notation is taken from \cite{G4}.
Therefore, to prove the discrete Laplacian $\Delta$ is negative semi-definite,
we just need to prove each matrix $(\Lambda_{\sigma})_{4\times 4}$ is negative semi-definite.
Following Glickenstein's approach in Appendix of \cite{G2}, we have the following result for the matrix $\Lambda_\sigma$.

\begin{proposition}\label{rank of Lambda-sigma}
Suppose $\sigma=\{ijkl\}$ is a nondegenerate   tetrahedron generated by Thurston's sphere packing metrics
with a weight $\Phi\equiv C\in [0, \frac{\pi}{2}]$.
Then $Rank (\Lambda_{\sigma})=3$ and the kernel of $\Lambda_{\sigma}$ is $\{t(1,1,1,1)|t\in \mathbb{R}\}$.
\end{proposition}

As the proof of Proposition \ref{rank of Lambda-sigma} is tedious direct calculations and not the main part of paper,
we put it in Appendix \ref{appendix}.\\

\noindent
\textbf{Proof of Theorem \ref{main theorem positivity of Laplacian} under the condition (\ref{condition 2})}
As $\Lambda=\sum_{\sigma\in T}(\Lambda_{\sigma})_{E}$, to prove $\Lambda$ is negative semi-definite,
we just need to prove that the three nonzero eigenvalues of $\Lambda_\sigma$ are negative for any $\sigma\in T$
and any $(\eta, r)\in \widetilde{\Omega}_{\sigma}$,
where
\begin{equation*}
\begin{aligned}
\widetilde{\Omega}_{\sigma}=\cup_{\eta\in \widetilde{W}}(\eta\times\Omega_{1234}(\eta))=\{(\eta, r)\in \widetilde{W}\times \mathbb{R}^4_{>0}|Q(\eta, r)>0\}
\end{aligned}
\end{equation*}
with
\begin{equation*}
\begin{aligned}
\widetilde{W}=\{\eta\in [0,1]^6\setminus S|\Phi\equiv C\in [0, \frac{\pi}{2}]\}.
\end{aligned}
\end{equation*}
It is straight forward that $\widetilde{W}$ is connected, which implies $\widetilde{\Omega}_{\sigma}$ is
connected by Corollary \ref{general connectness of 10-dim admissible space}.
Especially, for $\overline{\eta}=(0,0,0,0,0,0)$ and $\overline{r}=(1,1,1,1)$, we have $(\overline{\eta}, \overline{r})\in \widetilde{\Omega}_\sigma$.
By (\ref{calculation of h_ij,k}) and (\ref{calculation of h_ijk,l}), we have $A_{ij,kl}(\overline{\eta}, \overline{r})>0$, which implies $\Lambda_\sigma(\overline{\eta}, \overline{r})$
is negative semi-definite with rank $3$ by the the proof in Subsubsection \ref{subsection simple case of negative Lap}.
Therefore three nonzero eigenvalues of $\Lambda_\sigma(\overline{\eta}, \overline{r})$ are negative.

By Proposition \ref{rank of Lambda-sigma},
the rank of $\Lambda_\sigma(\eta, r)$ is 3 for any $(\eta, r)\in \widetilde{\Omega}_\sigma$.
Combining the continuity of eigenvalues and the connectedness of $\widetilde{\Omega}_\sigma$ in Corollary \ref{general connectness of 10-dim admissible space},
we have three nonzero eigenvalues of $\Lambda_\sigma$ are negative for any $(\eta, r)\in \widetilde{\Omega}_{\sigma}$.
This implies the negative semi-definiteness of $\Lambda_\sigma$ and $\Lambda=\sum_{\sigma\in T}(\Lambda_{\sigma})_{E}$.

The kernel of $\Lambda$ is $\{t(1,\cdots,1)|t\in \mathbb{R}\}$ follows from the connectedness of the manifold $M$ and
the kernel of $\Lambda_\sigma$ is $\{t(1,1,1,1)|t\in \mathbb{R}\}$ for any $\sigma\in T$. \qed

\section{Rigidity of Thurston's sphere packing}\label{section 4}

Recall the definition of discrete Hilbert-Einstein functional
$$F=\sum_{i\sim j}K_{ij}l_{ij}=\sum_{i\sim j}(2\pi-\sum \beta_{ij,kl})l_{ij}.$$
By the Schl\"{a}fli formula, we have
$dF=\sum_{i\sim j}K_{ij}dl_{ij}$,
which implies
$$\frac{\partial F}{\partial r_i}=\sum_{j\sim i}K_{ij}\frac{\partial l_{ij}}{\partial r_i}=K_i.$$
Further calculations show
\begin{equation*}
\begin{aligned}
\frac{\partial^2 F}{\partial r_i\partial r_j}=\sum_{s\sim t}\sum_{u\sim v} \frac{\partial l_{uv}}{\partial r_j}\cdot \frac{\partial K_{st}}{\partial l_{uv}}\cdot \frac{\partial l_{st}}{\partial r_i} +\sum_{s\sim t}K_{st}\frac{\partial l_{st}}{\partial r_i\partial r_j},
\end{aligned}
\end{equation*}
which implies
\begin{equation}\label{Hess F tempo}
\begin{aligned}
Hess_r F=-\sum_{\sigma=\{ijkl\}\in T}\left[R_\sigma D_{\sigma}^T(\frac{\partial \beta}{\partial l})_\sigma D_{\sigma} R_\sigma\right]_E
+\sum_{i\sim j} K_{ij}Hess_r(l_{ij}),
\end{aligned}
\end{equation}
where
$$R_\sigma=diag(r_i^{-1}, r_j^{-1}, r_k^{-1}, r_l^{-1}),$$
\begin{equation*}
\begin{aligned}
\left(\frac{\partial \beta}{\partial l}\right)_\sigma=
\left(
  \begin{array}{cccccc}
    \frac{\partial \beta_{ij}}{\partial l_{ij}} & \frac{\partial \beta_{ij}}{\partial l_{ik}} & \frac{\partial \beta_{ij}}{\partial l_{il}} & \frac{\partial \beta_{ij}}{\partial l_{jk}} & \frac{\partial \beta_{ij}}{\partial l_{jl}} & \frac{\partial \beta_{ij}}{\partial l_{kl}} \\
    \frac{\partial \beta_{ik}}{\partial l_{ij}} & \frac{\partial \beta_{ik}}{\partial l_{ik}} & \frac{\partial \beta_{ik}}{\partial l_{il}} & \frac{\partial \beta_{ik}}{\partial l_{jk}} & \frac{\partial \beta_{ik}}{\partial l_{jl}} & \frac{\partial \beta_{ik}}{\partial l_{kl}} \\
    \frac{\partial \beta_{il}}{\partial l_{ij}} & \frac{\partial \beta_{il}}{\partial l_{ik}} & \frac{\partial \beta_{il}}{\partial l_{il}} & \frac{\partial \beta_{il}}{\partial l_{jk}} & \frac{\partial \beta_{il}}{\partial l_{jl}} & \frac{\partial \beta_{il}}{\partial l_{kl}} \\
    \frac{\partial \beta_{jk}}{\partial l_{ij}} & \frac{\partial \beta_{jk}}{\partial l_{ik}} & \frac{\partial \beta_{jk}}{\partial l_{il}} & \frac{\partial \beta_{jk}}{\partial l_{jk}} & \frac{\partial \beta_{jk}}{\partial l_{jl}} & \frac{\partial \beta_{jk}}{\partial l_{kl}} \\
    \frac{\partial \beta_{jl}}{\partial l_{ij}} & \frac{\partial \beta_{jl}}{\partial l_{ik}} & \frac{\partial \beta_{jl}}{\partial l_{il}} & \frac{\partial \beta_{jl}}{\partial l_{jk}} & \frac{\partial \beta_{jl}}{\partial l_{jl}} & \frac{\partial \beta_{jl}}{\partial l_{kl}} \\
    \frac{\partial \beta_{kl}}{\partial l_{ij}} & \frac{\partial \beta_{kl}}{\partial l_{ik}} & \frac{\partial \beta_{kl}}{\partial l_{il}} & \frac{\partial \beta_{kl}}{\partial l_{jk}} & \frac{\partial \beta_{kl}}{\partial l_{jl}} & \frac{\partial \beta_{kl}}{\partial l_{kl}} \\
  \end{array}
\right),
D_\sigma=\left(
                                                                            \begin{array}{cccc}
                                                                              d_{ij} & d_{ji} & 0 & 0 \\
                                                                              d_{ik} & 0 & d_{ki} & 0 \\
                                                                              d_{il} & 0 & 0 & d_{li} \\
                                                                              0 & d_{jk} & d_{kj} & 0 \\
                                                                              0 & d_{jl} & 0 & d_{lj} \\
                                                                              0 & 0 & d_{kl} & d_{lk} \\
                                                                            \end{array}
                                                                          \right).
\end{aligned}
\end{equation*}
Recall the following important result of Glickenstein.

\begin{theorem}[\cite{G3}, Theorem 31]\label{Glickenstein's variational formula}
Suppose $\sigma=\{ijkl\}$ is a nondegenerate tetrahedron generated by a sphere packing metric $r: V\rightarrow (0,+\infty)$.
Set $u_i=\ln r_i$. Then we have
\begin{description}
  \item[(1)] if $j\sim i$, then
\begin{equation*}
\sum_{k\sim i}d_{ik}\frac{\partial \beta_{ik}}{\partial u_j}=\frac{2A_{ij, kl}}{l_{ij}},
\end{equation*}
  \item[(2)]
\begin{equation*}
\sum_{k\sim i}d_{ik}\frac{\partial \beta_{ik}}{\partial u_i}=-\sum_{j\sim i}\frac{2A_{ij, kl}}{l_{ij}}.
\end{equation*}
\end{description}
\end{theorem}

\begin{remark}
For the case of Euclidean tangential sphere packings on 3-dimensional manifolds,
Theorem \ref{Glickenstein's variational formula} was first proved by Glickenstein \cite{G1}.
Theorem \ref{Glickenstein's variational formula} is valid for much more general discrete conformal variations introduced by Glickenstein \cite{G3},
including the inversive distance sphere packing \cite{BS, G3, G4}
and vertex scaling \cite{L1,RW} on three dimensional manifolds. Please refer to \cite{G3, GT} for more details.
\end{remark}

Applying Theorem \ref{Glickenstein's variational formula} to (\ref{Hess F tempo}), we have
\begin{equation*}
\begin{aligned}
Hess_r F=-2\sum_{\sigma\{ijkl\}\in T}\left[R_\sigma \Lambda_\sigma R_\sigma\right]_E+\sum_{i\sim j} K_{ij}Hess_r(l_{ij})=-2R\Lambda R +\sum_{i\sim j} K_{ij}Hess_r(l_{ij}),
\end{aligned}
\end{equation*}
where $R=diag\{r_1^{-1}, \cdots, r_{|V|}^{-1}\}$ and $\Lambda$ is the matrix corresponding to the discrete Laplacian.
By direct calculations, we have
$$Hess_r(l_{ij})=\frac{\sin^2\Phi_{ij}}{l_{ij}^3}\left(
                                                   \begin{array}{cc}
                                                     r_j^2 & -r_ir_j \\
                                                     -r_ir_j & r_i^2 \\
                                                   \end{array}
                                                 \right)_E
$$
embedded as a $|V|\times |V|$ matrix, which is positive semi-definite and has $r=(r_1, \cdots, r_{|V|})$ as a null vector.
Therefore,
\begin{equation}\label{hess of F}
\begin{aligned}
Hess_r F=-2R\Lambda R +\sum_{i\sim j}\frac{K_{ij}\sin^2\Phi_{ij}}{l_{ij}^3}\left(
                                                   \begin{array}{cc}
                                                     r_j^2 & -r_ir_j \\
                                                     -r_ir_j & r_i^2 \\
                                                   \end{array}
                                                 \right)_E.
\end{aligned}
\end{equation}
By Theorem \ref{main theorem positivity of Laplacian},
$r=(r_1, \cdots, r_{|V|})$ is in the kernel of $Hess_r F$, which implies $Hess_r F$ has a zero eigenvalue for any sphere packing metric.\\

\noindent
\textbf{Proof of Theorem \ref{main theorem infinitesimal rigidity for scalar curvature}}~ ~
If $\overline{r}$ is a sphere packing metric with $K_{ij}(\overline{r})\sin^2 \Phi_{ij}\geq 0, \forall \{ij\}\in E$,
then
$Hess_r F(\overline{r})\geq 0$
and the kernel of $Hess_rF(\overline{r})$ is $\{t(\overline{r}_1, \cdots, \overline{r}_{|V|})|t\in \mathbb{R}\}$
by (\ref{hess of F}) and Theorem \ref{main theorem positivity of Laplacian}.
Therefore, $Hess_rF(\overline{r})$ has a zero eigenvalue and $|V|-1$ nonzero eigenvalues.
Note that $Hess_rF(r)$ has a zero eigenvalue with eigenvector $r=(r_1, \cdots, r_{|V|})$ for any sphere packing metric $r$.
Therefore, there exists a convex neighborhood $U$ of $\overline{r}$ such that
$Hess_r F(r)\geq 0$
and the kernel of $Hess_rF(r)$ is $\{t(r_1, \cdots, r_{|V|})|t\in \mathbb{R}\}$ for any $r\in U$.

Suppose there are two different sphere packing metrics $r_A$ and $r_B$ in $U$ such that $K_i(r_A)=K_i(r_B), \forall i\in V$.
Set
$f(t)=F(tr_A+(1-t)r_B).$
Then we have
$$f'(t)=\nabla F\cdot (r_A-r_B)=\sum_i (r_{A,i}-r_{B,i})K_i(tr_A+(1-t)r_B).$$
By $K_i(r_A)=K_i(r_B), \forall i\in V$, we have $f'(0)=f'(1)$.
Furthermore,
$$f''(t)=(r_A-r_B)^T\cdot Hess_r F|_{(tr_A+(1-t)r_B)}\cdot (r_A-r_B)\geq 0$$
by $Hess_r F\geq 0$ in $U$, which implies $f'(t)$ is nondecreasing for $t\in [0,1]$.
By $f'(0)=f'(1)$, we have $f'(t)=f'(0), \forall t\in [0,1]$ and $f''(t)\equiv 0$ in $[0,1]$,
which implies $r_A=cr_B$ by the fact that the kernel of $Hess_rF(r)$ is $\{t(r_1, \cdots, r_{|V|})|t\in \mathbb{R}\}$  for any $r\in U$.
\qed

\noindent


For a generic Thurston's sphere packing metric, we have the following result,
which is equivalent to Theorem \ref{main theorem infinitesimal rigidity}.

\begin{theorem}\label{main rigidity theorem}
Suppose $(M, \mathcal{T})$ is a triangulated closed connected 3-manifold
with a weight $\Phi: E\rightarrow [0,\frac{\pi}{2}]$ satisfying (\ref{condition 1}) or (\ref{condition 2}).
Any nondegenerate Thurston's   sphere packing metric $\overline{r}\in \Omega$ admits
a neighborhood $U\subset \Omega$ such that
if $r\in U$ has the same combinatorial Ricci curvature as $\overline{r}$, then $r=c\overline{r}$ for some positive constant $c$.
\end{theorem}

\proof
We take the modified discrete Hilbert-Einstein functional
$$\overline{F}=F-\sum K_{ij}(\overline{r})l_{ij}=\sum K_{ij}l_{ij}-\sum K_{ij}(\overline{r})l_{ij}$$
as a function of $r=(r_1, \cdots, r_N)$.
By direct calculations, we have
$$Hess_r\overline{F}
=-2R\Lambda R
+\sum (K_{ij}-K_{ij}(\overline{r}))\frac{\sin^2\Phi_{ij}}{l_{ij}^3}\left(
                                                   \begin{array}{cc}
                                                     r_j^2 & -r_ir_j \\
                                                     -r_ir_j & r_i^2 \\
                                                   \end{array}
                                                 \right)_E.$$
It is straight forward to check that $Hess_{r}\overline{F}\cdot r=0$, which implies that $Rank(Hess_{r}\overline{F})\leq |V|-1$.
Furthermore,
$Hess_r\overline{F}|_{\overline{r}}
=-2R\Lambda R$,
which is positive semi-definite with kernel $\{t\overline{r}|t\in \mathbb{R}\}$. This implies that $Rank(Hess_{r}\overline{F})=|V|-1$ at $\overline{r}$.

Note that the rank of a symmetric matrix is the number of nonzero eigenvalues.
By the continuity of the eigenvalues and $Rank(Hess_{r}\overline{F})\leq |V|-1$,
there is a convex neighborhood $U$ of $\overline{r}$ such that
$Rank(Hess_{r}(\overline{F}))=|V|-1$ and $Hess_{r}\overline{F}$ is positive semi-definite with kernel $\{tr|t\in \mathbb{R}\}$ for any $r\in U$.

Note that
$$\nabla_{r_i}\overline{F}=\sum_{j\sim i} (K_{ij}(r)-K_{ij}(\overline{r}))\frac{\partial l_{ij}}{\partial r_i}
=\sum_{j\sim i} (K_{ij}(r)-K_{ij}(\overline{r}))\cos \tau_{ij}.$$
If $K_{ij}(r)=K_{ij}(\overline{r})$ for some $r\in U$, then we have
$\nabla \overline{F}(r)=\nabla \overline{F}(\overline{r})=0.$
Define
$f(t)=\overline{F}(tr+(1-t)\overline{r})=\overline{F}(\overline{r}+t(r-\overline{r})),$
then $f(t)$ is a convex function on $[0,1]$ by $f''(t)=(r-\overline{r})^T\cdot Hess_r\overline{F}\cdot (r-\overline{r})\geq 0$.
$f'(t)=\nabla\overline{F}(tr+(1-t)\overline{r})\cdot (r-\overline{r})$
implies $f'(0)=f'(1)=0$, which
gives $f'(t)\equiv 0$ on $[0,1]$.
Then we have
$f''(t)=(r-\overline{r})^T\cdot Hess_r\overline{F}\cdot (r-\overline{r})=0,$
which implies that $r=c\overline{r}$ by $Hess_r\overline{F}$ is positive semi-definite
with kernel $\{tr|t\in \mathbb{R}\}$.\qed

%

\section{Appendix: Proof of Proposition \ref{rank of Lambda-sigma}}\label{appendix}
In this appendix, we describe a proof of Proposition \ref{rank of Lambda-sigma}.
By the definition of $\Lambda_\sigma$, $(1,1,1,1)$ is in the kernel of $\Lambda_\sigma$. We just need to
show that $Rank(\Lambda_\sigma)=3$. Set
\begin{equation*}
\begin{aligned}
I_{ijk}=&r_i^2r_j^2(1-\eta_{ij}^2)+r_i^2r_jr_k\Gamma_{ijk}
+r_ir_j^2r_k\Gamma_{jik},\\
B_{ijkl}=&2l_{ij}^2(l_{ik}^2+l_{il}^2-l_{kl}^2)
-(l_{ij}^2+l_{ik}^2-l_{jk}^2)(l_{ij}^2+l_{il}^2-l_{jl}^2),
\end{aligned}
\end{equation*}
then (\ref{calculation of h_ij,k}) and (\ref{cos beta_ij,kl}) can be rewritten as
\begin{equation}\label{hijknew}
\begin{aligned}
h_{ij,k}=\frac{I_{ijk}}{2l_{ij}A_{ijk}},
\cos\beta_{ij}=\frac{B_{ijkl}}{16A_{ijk}A_{ijl}}.
\end{aligned}
\end{equation}
Moreover, $h_{ijk,l}$ can be calculated as
\begin{eqnarray*}
h_{ijk,l}=\frac{h_{ij,l}-h_{ij,k}\cos\beta_{ij}}{\sin\beta_{ij}}
=\frac{16I_{ijl}A_{ijk}^2-I_{ijk}B_{ijkl}}{48l_{ij}^2A_{ijk}V_{ijkl}},
\end{eqnarray*}
where we have used the fact
$\sin\beta_{ij}=\frac{3V_{ijkl}l_{ij}}{2A_{ijk}A_{ijl}}.$
Set
$C_{ijkl}=16I_{ijl}A_{ijk}^2-I_{ijk}B_{ijkl},$
then we have
\begin{eqnarray}
h_{ijk,l}=\frac{C_{ijkl}}{48l_{ij}^2A_{ijk}V_{ijkl}}.\label{hijklnew}
\end{eqnarray}
Combining (\ref{hijknew}), (\ref{hijklnew}) with(\ref{Aijklnew}) gives
\begin{eqnarray*}
A_{ij,kl}=\frac{I_{ijk}C_{ijkl}A_{ijl}^2+I_{ijl}C_{ijlk}A_{ijk}^2}
{192l_{ij}^3A_{ijk}^2A_{ijl}^2V_{ijkl}}.
\end{eqnarray*}
Let
$U_{\sigma}=-2\Lambda_{\sigma}V_{ijkl}.$
Then
\begin{eqnarray*}
(U_{\sigma})_{ij}=-\frac{I_{ijk}C_{ijkl}A_{ijl}^2+I_{ijl}C_{ijlk}A_{ijk}^2}
{96l_{ij}^4A_{ijk}^2A_{ijl}^2}, \ i\sim j,
\end{eqnarray*}
the numerator and denominator of which are polynomials of $r_i,r_j,r_k,r_l$.

In order to show $Rank(\Lambda_\sigma)=3$,
we just need to prove that a $3$-dimensional minor of $U_\sigma$ is nonzero, which is
a complicated symbolic calculation.
Denote the determinant of the $(1, 2)$ minor of the matrix $U_\sigma$ by $M_{12}$,
which is the submatrix of $U_\sigma$ with the $1$-th row and the $2$-th column removed.
We will show that $M_{12}<0$ if the Euclidean tetrahedron is nondegenerate as follows.

Let
$x=\cos\Phi=\cos C,$
then $x\in[0,1]$ for $C\in[0,\frac{\pi}{2}]$.
Set
\begin{eqnarray*}
\widetilde{M}_{12}=884736A_{ijk}^4A_{ijl}^4A_{ikl}^4A_{jkl}^4l_{ij}^4l_{ik}^4
l_{il}^4l_{jk}^4l_{jl}^4l_{kl}^4 M_{12}
\end{eqnarray*}
to kill the denominator of $M_{12}$.
Without loss of generality, one can set $r_i=1, r_j=r_2, r_k=r_3, r_l=r_4$ by the scaling property of $\Lambda_\sigma$ and $U_\sigma$.
Direct calculations show
that the nondegenerate condition $\det G_0>0$ is
$\det G_0=8(1+x)^2M_1>0$,
where
\begin{equation*}
\begin{aligned}
M_1=&(1-2x)\bigg{[}r_2^2r_3^2+r_2^2r_4^2+r_3^2r_4^2+r_2^2r_3^2r_4^2\bigg{]}\\
&+2x\bigg{[}r_2^2r_3r_4+r_2r_3^2r_4+r_2^2r_3^2r_4
+r_2r_3r_4^2+r_2^2r_3r_4^2+r_2r_3^2r_4^2\bigg{]}.
\end{aligned}
\end{equation*}
Therefore, the nondegenerate condition is equivalent to $M_1>0$.
The calculations of $\widetilde{M}_{12}$ and $M_1$ are accomplished with the help of Mathematica.
We provide the full  Mathematica codes and the relevant outputs of those codes as an ancillary file (\texttt{determinant1.pdf})
with the arXiv version of this paper  \cite{HX0}.

Set
\begin{eqnarray*}
\widetilde{M}_{12}=q_{12}\cdot
\frac{l_{ij}^2l_{ik}^2l_{il}^2l_{jk}^2l_{jl}^2l_{kl}^2r_i^2r_j^2r_k^2r_l^2
(1+x)^9}{16384}+R_{12}.
\end{eqnarray*}
Direct calculations with the help of Mathematica show $R_{12}=0$, which implies  $M_{12}<0$ is equivalent
to $q_{12}<0$ under the condition that the tetrahedron is nondegenerate.
The decomposition of $q_{12}$
is accomplished again with the help of Mathematica. The result is
\begin{eqnarray*}
q_{12}=M_1^2\sum_{i,j,k=0}^{14}C_{i,j,k}(x)r_2^ir_3^jr_4^k.
\end{eqnarray*}
With the help of Mathematica,
we give all the graphs and expressions of the coefficient functions $C_{i,j,k}(x)(i,j,k=0,\cdots 14)$ for $x\in [0,1]$,
which shows that the coefficients $C_{i,j,k}(x)(i,j,k=0,\cdots 14)$ are non-positive when $x\in [0,1]$ and not all zero at the same value of $x$.
Therefore, we have $q_{12}<0$ if the sphere packing metric is nondegenerate, which completes the proof of Proposition \ref{rank of Lambda-sigma}.
For the calculations of $R_{12}$, decomposition of $q_{12}$ and graphs and expressions of $C_{i,j,k}(x)$,
we provide the full  Mathematica codes and the relevant outputs of those codes as an ancillary file (\texttt{determinant2.pdf}) with
the arXiv version of this paper \cite{HX0}.
\qed
\bigskip

Here we give the full expression of $q_{12}$ for some special weights.
When $\Phi\equiv 0$,
\begin{equation*}
\begin{aligned}
q_{12}=-6291456 M_1^2 r_2^5(1+r_2)^2r_3^5(1+r_3)^2(r_2+r_3)^2(1+r_2+r_3)r_4^5
(1+r_4)^2\nonumber\\
\cdot (r_2+r_4)^2(1+r_2+r_4)(r_3+r_4)^2(1+r_3+r_4)(r_2+r_3+r_4),
\end{aligned}
\end{equation*}
which coincides with the result obtained by Glickenstein \cite{G1, G2} for tangential sphere packing metrics.
If $\Phi\equiv \frac{\pi}{3}$,
\begin{equation*}
\begin{aligned}
q_{12}=&-\frac{243}{8}M_1^2(r_2+r_3+r_2r_3)^4(r_2+r_4+r_2r_4)^2(r_3+r_4+r_3r_4)^4
(r_3r_4+r_2r_3+r_2r_4)^4\\
&\cdot [r_3r_4(1+r_3+r_4)+r_2^2(r_3+r_4+r_3r_4)+
r_2(r_3+r_3^2+r_4+r_3^2r_4+r_4^2+r_3r_4^2)].
\end{aligned}
\end{equation*}

In Proposition \ref{rank of Lambda-sigma},
we only prove the case that the intersection angles $\Phi_{ij}$ are all the same for any edge $\{ij\}\in E$
due to the computer performance.
However, it is conceivable that Proposition \ref{rank of Lambda-sigma} is true for any weight $\Phi: E\rightarrow [0, \frac{\pi}{2}]$.
If Proposition \ref{rank of Lambda-sigma} is true for any weight $\Phi: E\rightarrow [0,\frac{\pi}{2}]$,
then Theorem \ref{main theorem infinitesimal rigidity} and Theorem \ref{main theorem infinitesimal rigidity for scalar curvature}
 are true for any weight $\Phi: E\rightarrow [0,\frac{\pi}{2}]$.

(Xiaokai He) School of Mathematics and Statistics, Hunan First Normal University, Changsha 410205, P.R. China

E-mail:hexiaokai77@163.com \\[2pt]

(Xu Xu) School of Mathematics and Statistics, Wuhan University, Wuhan 430072, P.R. China

E-mail: xuxu2@whu.edu.cn\\[2pt]

\end{document}